# ESTIMATES OF MOMENTS AND TAILS OF GAUSSIAN CHAOSES[1]

By Rafał Latała

*Warsaw University*

We derive two-sided estimates on moments and tails of Gaussian chaoses, that is, random variables of the form $\sum a_{i_1,\ldots,i_d} g_{i_1} \cdots g_{i_d}$, where $g_i$ are i.i.d. $\mathcal{N}(0,1)$ r.v.'s. Estimates are exact up to constants depending on $d$ only.

**1. Introduction.** The purpose of this paper is to give precise bounds on moments and tails of Gaussian chaoses of order $d$, that is, random variables of the form $S = \sum_{i_1 < i_2 < \cdots < i_d} a_{i_1,\ldots,i_d} g_{i_1} \cdots g_{i_d}$. In the sequel, we will only consider decoupled chaoses $\tilde{S} = \sum_{\mathbf{i}} a_{\mathbf{i}} g_{i_1}^{(1)} \cdots g_{i_d}^{(d)}$, where $g_i^{(k)}$ are independent standard $\mathcal{N}(0,1)$ normal random variables and $(a_{\mathbf{i}}) = (a_{i_1,\ldots,i_d})_{1 \le i_1,\ldots,i_d \le N}$ is a finite multi-indexed matrix—under natural symmetry assumptions, moments and tails of $S$ and $\tilde{S}$ are comparable with constants depending only on $d$ (cf. [5]).

For $d = 1$, we obviously have for $p \geq 2$,

$$(1) \qquad \left\| \sum_i a_i g_i \right\|_p = \left( \sum_i a_i^2 \right)^{1/2} \|g\|_p \sim \sqrt{p} \left( \sum_i a_i^2 \right)^{1/2}.$$

For the chaoses of order 2, we have for any finite rectangular matrix $(a_{ij})$ and $p \geq 2$,

$$(2) \qquad \left\| \sum_{ij} a_{ij} g_i^{(1)} g_j^{(2)} \right\|_p \sim \sqrt{p} \|(a_{ij})\|_{\{1,2\}} + p \|(a_{ij})\|_{\{1\}\{2\}},$$

where $\|(a_{ij})\|_{\{1,2\}} := \|(a_{ij})\|_{\mathrm{HS}} = (\sum_{ij} a_{ij}^2)^{1/2}$ and

$$\|(a_{ij})\|_{\{1\}\{2\}} := \sup\left\{ \sum_{ij} a_{ij} x_i y_j : \|x\|_2 \le 1, \|y\|_2 \le 1 \right\}.$$

Received May 2005; revised December 2005.
[1]Supported in part by KBN Grant 2 PO3A 027 22.

The upper part of the estimate (2) was obtained in [6]; the lower one is much easier (cf. [7]). One of the reasons why the case where $d=2$ turned out to be relatively simple is that every square matrix is orthogonally equivalent to the diagonal matrix.

For $d \geq 3$, Borell [4] and Arcones and Giné [3] showed that

$$(3) \quad \|S\|_p \stackrel{d}{\sim} \sum_{k=1}^{d} p^{k/2} \mathbf{E} \sup \left\{ \sum_{\mathbf{i}} a_{\mathbf{i}} \prod_{l=1}^{k} x_{i_k}^{(k)} \prod_{l=k+1}^{d} g_{i_l}^{(l)} : \|x^{(l)}\|_2 \leq 1, 1 \leq l \leq k \right\}.$$

The above formula gives the precise dependence on $p$, but unfortunately involves suprema of empirical processes that are, in general, not easy to estimate. (For generalizations of the above formula to the non-Gaussian case, cf. [1] and [10].) In this paper, we present bounds on moments and tails that involve only deterministic quantities.

The paper is organized as follows. In the next section, we present notation and definitions that will be used in the rest of the paper and formulate main results. In Section 3 we obtain bounds on entropy numbers for distances on products of Euclidean balls. This will provide a crucial tool to estimate suprema of certain Gaussian processes that naturally appear in the study of Gaussian chaoses. Finally, in the last section, we present proofs of main results.

**2. Notation and main results.** We use the letter $C$ to denote universal positive constants that may change from occurrence to occurrence and $C(d)$ to denote positive constants depending only on $d$. [$C(d)$ may also differ at each occurrence.] We write $f \sim g$ if $\frac{1}{C} f \leq g \leq Cf$ and $f \sim_d g$ if $\frac{1}{C(d)} f \leq g \leq C(d)f$. The canonical Euclidean norm of a vector $x$ is denoted by $\|x\|_2$. Recall that the $p$th norm of a real random variable $X$ is defined as $\|X\|_p := (\mathbf{E}|X|^p)^{1/p}$.

Let $d \geq 1$ and $A = (a_{\mathbf{i}})_{1 \leq i_1, \ldots, i_d \leq n}$ be a finite multi-indexed matrix of order $d$. If $\mathbf{i} \in \{1, \ldots, n\}^d$ and $I \subset \{1, \ldots, d\}$, then we define $i_I := (i_j)_{j \in I}$. For disjoint nonempty subsets $I_1, \ldots, I_k$ of $\{1, \ldots, d\}$, we put

$$\|A\|_{I_1, \ldots, I_k} := \sup \left\{ \sum_{\mathbf{i}} a_{\mathbf{i}} x_{i_{I_1}}^{(1)} \cdots x_{i_{I_k}}^{(k)} : \sum_{i_{I_1}} (x_{i_{I_1}}^{(1)})^2 \leq 1, \ldots, \sum_{i_{I_k}} (x_{i_{I_k}}^{(k)})^2 \leq 1 \right\}.$$

Thus, for example,

$$\|A\|_{\{1,\ldots,d\}} = \left( \sum_{\mathbf{i}} a_{\mathbf{i}}^2 \right)^{1/2}$$

and

$$\|(a_{ijk})\|_{\{1\}\{2,3\}} = \sup \left\{ \left( \sum_{jk} \left( \sum_i a_{ijk} x_i \right)^2 \right)^{1/2} : \sum_i x_i^2 \leq 1 \right\}$$

$$= \sup\left\{\left(\sum_i\left(\sum_{jk}a_{ijk}x_{jk}\right)^2\right)^{1/2} : \sum_{jk}x_{jk}^2 \leq 1\right\}.$$

By $S(k,d)$, we denote a set of all partitions of $\{1,\ldots,d\}$ into $k$ nonempty disjoint sets $I_1,\ldots,I_k$. For $p \geq 1$, we put

$$m_p(A) := \sum_{k=1}^d p^{k/2} \sum_{(I_1,\ldots,I_k)\in S(k,d)} \|A\|_{I_1,\ldots,I_k}.$$

Our main result is the following:

THEOREM 1. *For any multi-indexed finite matrix* $A = (a_{\mathbf{i}})_{1\leq i_1,\ldots,i_d\leq n}$ *and* $p \geq 2$, *we have*

(4) $$\frac{1}{C(d)}m_p(A) \leq \left\|\sum_{\mathbf{i}} a_{\mathbf{i}} \prod_{j=1}^d g_{i_j}^{(j)}\right\|_p \leq C(d)m_p(A).$$

Theorem 1 may be easily translated into the following two-sided estimate for tails:

COROLLARY 1. *For any* $t > 0$ *and* $d \geq 2$, *we have*

$$\frac{1}{C(d)}\exp\left[-C(d)\min_{1\leq k\leq d}\min_{(I_1,\ldots,I_k)\in S(k,d)}\left(\frac{t}{\|A\|_{I_1,\ldots,I_k}}\right)^{2/k}\right]$$
$$\leq \mathbf{P}\left(\left|\sum_{\mathbf{i}} a_{\mathbf{i}} \prod_{j=1}^d g_{i_j}^{(j)}\right| \geq t\right)$$
$$\leq C(d)\exp\left[-\frac{1}{C(d)}\min_{1\leq k\leq d}\min_{(I_1,\ldots,I_k)\in S(k,d)}\left(\frac{t}{\|A\|_{I_1,\ldots,I_k}}\right)^{2/k}\right].$$

In view of (3), it is clear that the proof of (4) should be based on the estimation of norms of some random Gaussian matrices. The next theorem is, in our opinion, of independent interest and has recently been applied in [2] to obtain moment estimates for canonical $U$-statistics.

THEOREM 2. *For any* $d \geq 2$ *and any finite matrix* $A$, *we have for* $p \geq 2$,

(5) $$\mathbf{E}\left\|\left(\sum_{i_d} a_{\mathbf{i}} g_{i_d}\right)\right\|_{\{1\}\ldots\{d-1\}} \leq C(d)p^{(1-d)/2}m_p(A).$$

REMARK. We suspect that a stronger estimate may actually hold, namely

$$\mathbf{E}\left\|\left(\sum_{\mathbf{i}_d} a_{\mathbf{i}} g_{i_d}\right)\right\|_{\{1\}\ldots\{d-1\}}$$
$$\leq C(d) \inf_{p \geq 1} (s_{d-1}(A) + p^{1/2} \|A\|_{\{1\}\ldots\{d\}} + p^{(2-d)/2} \|A\|_{\{1,\ldots,d\}})$$
$$\leq C(d)(s_{d-1}(A) + \|A\|_{\{1\}\ldots\{d\}}^{(d-2)/(d-1)} \|A\|_{\{1,\ldots,d\}}^{1/(d-1)}),$$

where

$$s_{d-1}(A) := \sum_{1 \leq j \leq d-1} \|A\|_{\{j,d\},\{\{l\}\,:\,1 \leq l \leq d-1, l \neq j\}}.$$

However, we are not able to show this result for $d > 3$.

**3. Entropy estimates and Gaussian processes.** By $\gamma_{n,t}$, we will denote the distribution of $tG_n$, where $G_n = (g_1, \ldots, g_n)$ is a canonical $n$-dimensional Gaussian vector. We also put $G_n^{(i)} := (g_1^{(i)}, \ldots, g_n^{(i)})$ for i.i.d. copies of $G_n$.

If $\rho$ is a metric on a set $T$, $N(T, \rho, t)$ is the minimal number of closed balls of radius $t$ that are necessary to cover $T$. The closed unit Euclidean ball in $\mathbb{R}^n$ is denoted by $B_2^n$.

LEMMA 1. *For any norms $\alpha_1, \alpha_2$ on $\mathbb{R}^n$, $y \in S \subset B_2^n$ and $t > 0$,*

$$\gamma_{n,t}(x : \alpha_1(x-y) \leq 4t\mathbf{E}\alpha_1(G_n), \alpha_2(x) \leq 4t\mathbf{E}\alpha_2(G_n) + \alpha_2(y)) \geq \tfrac{1}{2} e^{-t^{-2}/2}.$$

PROOF. Let

$$K := \{x \in \mathbb{R}^n : \alpha_1(x) \leq 4t\mathbf{E}\alpha_1(G_n), \alpha_2(x) \leq 4t\mathbf{E}\alpha_2(G_n)\}.$$

By Chebyshev's inequality,

$$1 - \gamma_{n,t}(K) \leq \mathbf{P}(\alpha_1(tG_n) > 4\mathbf{E}\alpha_1(tG_n)) + \mathbf{P}(\alpha_2(tG_n) > 4\mathbf{E}\alpha_2(tG_n)) \leq 1/2.$$

By the symmetry of $K$, we obtain for any $y \in B_2^n$,

$$\gamma_{n,t}(y+K) = e^{-|y|^2/(2t^2)} \int_K e^{\langle y,x \rangle/t^2} \, d\gamma_{n,t}(x)$$
$$= e^{-|y|^2/(2t^2)} \int_K \tfrac{1}{2}(e^{\langle y,x \rangle/t^2} + e^{-\langle y,x \rangle/t^2}) \, d\gamma_{n,t}(x)$$
$$\geq \exp(-t^{-2}/2) \gamma_{n,t}(K) \geq \tfrac{1}{2} \exp(-t^{-2}/2).$$

Finally, note that if $x \in y + K$, then $\alpha_1(x-y) \leq 4t\mathbf{E}\alpha_1(G_n)$ and $\alpha_2(x) \leq \alpha_2(x-y) + \alpha_2(y) \leq 4t\mathbf{E}\alpha_2(G_n) + \alpha_2(y)$. □

Let $\alpha$ be a norm on $\mathbb{R}^{n_1\cdots n_d}$ and the distance $\rho_\alpha$ on $\mathbb{R}^{n_1}\times\cdots\times\mathbb{R}^{n_d}$ be defined by

$$\rho_\alpha(\mathbf{x},\mathbf{y}):=\alpha\left(\bigotimes_{i=1}^d x^i - \bigotimes_{i=1}^d y^i\right) \quad \text{for } \mathbf{x}=(x^1,\ldots,x^d), \mathbf{y}=(y^1,\ldots,y^d),$$

where $\bigotimes_{i=1}^d x^i := (\prod_{k=1}^d x^k_{i_k})_{i_1\leq n_1,\ldots,i_d\leq n_d}$. For $t>0$, $T\subset \mathbb{R}^{n_1}\times\cdots\times\mathbb{R}^{n_d}$, we put

$$W_d^T(\alpha,t):=\sum_{k=1}^d t^k \sum_{I\subset\{1,\ldots,d\}:\, \#I=k} W_I^T(\alpha),$$

where

$$W_I^T(\alpha):=\sup_{\mathbf{x}\in T}\mathbf{E}\alpha\left(\left(\prod_{k\notin I} x^k_{i_k}\prod_{k\in I} g^{(k)}_{i_k}\right)_{i_1,\ldots,i_d}\right).$$

To simplify the notation, we will write $W_d$ and $W_I$ instead of $W_d^T$ and $W_I^T$ if $T=B_2^{n_1}\times\cdots\times B_2^{n_d}$.

LEMMA 2. *For any $t>0$ and $\mathbf{x}\in B_2^{n_1}\times\cdots\times B_2^{n_d}$, we have*

(6) $$\gamma_{n_1+\cdots+n_d,t}(B_\alpha(\mathbf{x},W_d^{\{\mathbf{x}\}}(\alpha,4t)))\geq 2^{-d}\exp(-dt^{-2}/2).$$

PROOF. We will proceed by induction on $d$. For $d=1$, we have

$$B_\alpha(\mathbf{x},W_d^{\{\mathbf{x}\}}(\alpha,4t))=\{y\in\mathbb{R}^{n_1}:\alpha(x-y)\leq 4t\mathbf{E}\alpha(G_{n_1})\}$$

and (6) then follows by Lemma 1.

Now, suppose that (6) holds for $d-1$. We will show that it is also satisfied for $d$.

Let us first observe that

(7) $$\alpha\left(\bigotimes_{i=1}^d x^i - \bigotimes_{i=1}^d y^i\right)\leq \alpha^1(x^d-y^d)+\alpha_{y^d}\left(\bigotimes_{i=1}^{d-1} x^i - \bigotimes_{i=1}^{d-1} y^i\right),$$

where $\alpha^1$ and $\alpha_y$ are norms on $\mathbb{R}^{n_d}$ and $\mathbb{R}^{n_1\cdots n_{d-1}}$, respectively, defined by

$$\alpha^1(z):=\alpha\left(\bigotimes_{i=1}^{d-1} x^i\otimes z\right) \quad \text{and} \quad \alpha_y(z):=\alpha(z\otimes y).$$

Then, obviously,

(8) $$\mathbf{E}\alpha^1(G_{n_d})=W_{\{d\}}^{\{\mathbf{x}\}}(\alpha).$$

Moreover, if we put $\pi(\mathbf{x}) = (x^1, \ldots, x^{d-1})$ and define a norm $\alpha_t^2$ on $\mathbb{R}^{n_d}$ by the formula
$$\alpha_t^2(y) := W_{d-1}^{\{\pi(\mathbf{x})\}}(\alpha_y, t),$$
then

(9) $\qquad t\mathbf{E}\alpha_t^2(G_{n_d}) + \alpha_t^2(x^d) = \sum_{I \subset \{1,\ldots,d\} : I \neq \varnothing, \{d\}} t^{\#I} W_I^{\{\mathbf{x}\}}(\alpha, t).$

Observe also that by the induction assumption, we have for any $z \in \mathbb{R}^{n_d}$,

(10) $\qquad \gamma_{n_1+\cdots+n_{d-1},t}\left(\mathbf{y} \in \mathbb{R}^{n_1+\cdots+n_{d-1}} : \alpha_z\left(\bigotimes_{i=1}^{d-1} x^i - \bigotimes_{i=1}^{d-1} y^i\right) \leq \alpha_{4t}^2(z)\right)$
$$\geq 2^{1-d} \exp(-(d-1)t^{-2}/2).$$

Finally, let
$$A(\mathbf{x}) := \left\{\mathbf{y} \in \mathbb{R}^{n_1+\cdots+n_d} : \alpha^1(x^d - y^d) \leq 4t\mathbf{E}\alpha^1(G_{n_d}),\right.$$
$$\alpha_{4t}^2(y^d) \leq 4t\mathbf{E}\alpha_{4t}^2(G_{n_d}) + \alpha_{4t}^2(x^d),$$
$$\left.\alpha_{y^d}\left(\bigotimes_{i=1}^{d-1} x^i - \bigotimes_{i=1}^{d-1} y^i\right) \leq \alpha_{4t}^2(y^d)\right\}.$$

By (7)–(9), we get $A(\mathbf{x}) \subset B_\alpha(\mathbf{x}, W_d^{\{\mathbf{x}\}}(\alpha, 4t))$ and, therefore, by (10), Lemma 1 and Fubini's theorem, we get
$$\gamma_{n_1+\cdots+n_d,t}(B_\alpha(\mathbf{x}, W_d^{\{\mathbf{x}\}}(\alpha, 4t))) \geq \gamma_{n_1+\cdots+n_d,t}(A(\mathbf{x})) \geq 2^{-d}\exp(-dt^{-2}/2). \square$$

COROLLARY 2. *For any $T \subset B_2^{n_1} \times \cdots \times B_2^{n_d}$ and $t \in (0,1]$,*
$$N(T, \rho_\alpha, W_d^T(\alpha, t)) \leq \exp(Cdt^{-2}).$$

*In particular,*
$$N(B_2^{n_1} \times \cdots \times B_2^{n_d}, \rho_\alpha, W_d(\alpha, t)) \leq \exp(Cdt^{-2}).$$

PROOF. Obviously, $W_d^T(\alpha, t) \geq \sup_{\mathbf{x} \in T} W_d^{\{\mathbf{x}\}}(\alpha, t)$. Therefore, by Lemma 2, we have for any $\mathbf{x} \in T$,

(11) $\qquad \gamma_{n_1+\cdots+n_d,t}(B_\alpha(\mathbf{x}, W_d^T(\alpha, 4t))) \geq 2^{-d}\exp(-dt^{-2}/2).$

Suppose that there exist $\mathbf{x}_1, \ldots, \mathbf{x}_N \in T$ such that $\rho_\alpha(\mathbf{x}_i, \mathbf{x}_j) > W_d^T(\alpha, t) \geq 2W_d^T(\alpha, t/2)$ for $i \neq j$. Then the sets $B_\alpha(\mathbf{x}_i, W_d^T(\alpha, t/2))$ are disjoint, so by (11), we obtain $N \leq 2^d \exp(32dt^{-2})$. Hence,
$$N(T, \rho_\alpha, W_d^T(\alpha, t)) \leq 2^d \exp(32dt^{-2}) \leq \exp(33dt^{-2}).$$

□

To finish this section, let us recall standard estimates for Gaussian processes.

LEMMA 3. *Let $(X_t)_{t\in T}$ be a centered Gaussian process and $T = \bigcup_{l=1}^{m} T_l$. Then*

$$\mathbf{E}\sup_{t\in T} X_t \leq \max_l \mathbf{E}\sup_{t\in T_l} X_t + C\sqrt{\log m}\sup_{t,s\in T}(\mathbf{E}(X_t - X_s)^2)^{1/2}.$$

PROOF. Obviously, $\mathbf{E}\sup_{t\in T} X_t = \mathbf{E}\max_l \sup_{t\in T_l}(X_t - X_{t_0})$ for any $t_0 \in T$. The lemma follows by integration by parts and the classical estimate (cf. [8], Theorem 7.1)

$$\mathbf{P}\left(\sup_{t\in T_l}(X_t - X_{t_0}) \geq \mathbf{E}\sup_{t\in T_l} X_t + u\sup_{t\in T_l}(\mathbf{E}(X_t - X_{t_0})^2)^{1/2}\right) \leq \exp(-u^2/2)$$

for $u > 0$. □

LEMMA 4. *Let $(X_t)_{t\in T}$ be a centered Gaussian process. Then for any $p \geq 2$,*

$$(12) \quad \frac{1}{C}\left(\left\|\sup_{t\in T} X_t\right\|_1 + \sqrt{p}\sigma\right) \leq \left\|\sup_{t\in T} X_t\right\|_p \leq \left\|\sup_{t\in T} X_t\right\|_1 + C\sqrt{p}\sigma,$$

*where $\sigma := \sup_{t\in T}(\mathbf{E}X_t^2)^{1/2}$.*

PROOF. The lower bound follows from the easy estimates

$$\left\|\sup_{t\in T} X_t\right\|_p \geq \left\|\sup_{t\in T}\max(X_t, 0)\right\|_p \geq \sup_{t\in T}\|\max(X_t, 0)\|_p \geq \sup_{t\in T}\|X_t\|_p/2$$

and the upper one by the concentration of suprema of Gaussian processes (cf. [8], Theorem 7.1) and integration by parts. □

**4. Proofs.** Let us start with some additional notation. For a matrix $A = (a_\mathbf{i})_{1\leq i_1,\ldots,i_d\leq n}$ of order $d \geq 2$, we set

$$s_{d-1}(A) := \sum_{1\leq j\leq d-1} \|A\|_{\{j,d\},\{\{l\}\,:\,1\leq l\leq d-1, l\neq j\}}$$

and for $1 \leq k \leq d-2$,

$$s_k(A) := \sum_{(I_1,\ldots,I_k)\in S(k,d)} \|A\|_{I_1,\ldots,I_k}.$$

On $\mathbb{R}^{(d-1)n} = (\mathbb{R}^n)^{d-1}$, we introduce the distance $\rho_A$ by the formula

$$\rho_A(\mathbf{x}, \mathbf{y}) := \left(\sum_{i_d}\left(\sum_{i_1,\ldots,i_{d-1}} a_{\mathbf{i}}\left(\prod_{k=1}^{d-1} x_{i_k}^k - \prod_{k=1}^{d-1} y_{i_k}^k\right)\right)^2\right)^{1/2},$$

where $\mathbf{x} = (x^1, \ldots, x^{d-1})$ and $\mathbf{y} = (y^1, \ldots, y^{d-1})$. We have

(13) $$\rho_A(\mathbf{x}, \mathbf{y}) = (\mathbf{E}(X_{\mathbf{x}} - X_{\mathbf{y}})^2)^{1/2},$$

where $X_{\mathbf{x}} := \sum_{i_1,\ldots,i_d} a_{\mathbf{i}} \prod_{k=1}^{d-1} x_{i_k}^k g_{i_d}$.

For $T \subset \mathbb{R}^{(d-1)n}$, we put

$$\Delta_A(T) := \sup\{\rho_A(\mathbf{x}, \mathbf{y}) : \mathbf{x}, \mathbf{y} \in T\}.$$

Let us note that in particular, we have

(14) $$\Delta_A((B_2^n)^{d-1}) \leq 2\sup\{\rho_A(\mathbf{x}, 0) : \mathbf{x} \in (B_2^n)^{d-1}\} = 2\|A\|_{\{1\}\ldots\{d\}}.$$

For a set $T \subset \mathbb{R}^{(d-1)n}$ and $I \subset \{1, \ldots, d-1\}$, we put

$$W_I^T(A) := \sup_{\mathbf{x} \in T}\left(\sum_{i_{I \cup \{d\}}}\left(\sum_{i_{\{1,\ldots,d-1\}\setminus I}} a_{\mathbf{i}} \prod_{k \in \{1,\ldots,d-1\}\setminus I} x_{i_k}^k\right)^2\right)^{1/2}$$

and for $1 \leq k \leq d-1$,

$$W_k^T(A) := \sum_{I \subset \{1,\ldots,d-1\}, \#I = k} W_I^T(A).$$

The next lemma shows how the results of the previous section may be adapted to the case of the particular metric $\rho_A$.

LEMMA 5.  *For any $0 < t \leq 1$ and $T \subset (B_2^n)^{d-1}$, we have*

(15) $$N\left(T, \rho_A, \sum_{k=1}^{d-1} t^k W_k^T(A)\right) \leq \exp(Cdt^{-2}).$$

*In particular,*

(16) $$N\left(T, \rho_A, tW_1^T(A) + \sum_{k=2}^{d-1} t^k s_{d-k}(A)\right) \leq \exp(Cdt^{-2}).$$

PROOF. Note that $\rho_A = \rho_\alpha$, where for $z \in \mathbb{R}^{n^{d-1}}$,

$$\alpha(z) := \left(\sum_{i_d}\left(\sum_{i_1,\ldots,i_{d-1}} a_{\mathbf{i}} z_{i_1,\ldots,i_{d-1}}\right)^2\right)^{1/2}.$$

We have for any $\mathbf{x} \in (\mathbb{R}^n)^{d-1}$ and $I \subset \{1, \ldots, d-1\}$,

$$\mathbf{E}\alpha\left(\prod_{k \in \{1,\ldots,d-1\}\setminus I} x_{i_k}^k \prod_{k \in I} g_{i_k}^{(k)}\right) \leq \left(\mathbf{E}\alpha^2\left(\prod_{k \in \{1,\ldots,d-1\}\setminus I} x_{i_k}^k \prod_{k \in I} g_{i_k}^{(k)}\right)\right)^{1/2}$$

$$= \left(\sum_{i_{I \cup \{d\}}} \left(\sum_{i_{\{1,\ldots,d-1\}\setminus I}} a_{\mathbf{i}} \prod_{k \in \{1,\ldots,d-1\}\setminus I} x_{i_k}^k\right)^2\right)^{1/2}.$$

Hence,

$$W_I^T(\alpha) \leq W_I^T(A)$$

and (15) immediately follows by Corollary 2.

Inequality (15) implies (16), since

$$W_I^T(A) \leq W_I^{(B_2^n)^{d-1}}(A) = \|A\|_{I \cup \{d\}, \{\{l\} : l \in \{1,\ldots,d-1\}\setminus I\}}. \qquad \square$$

We are now ready to present a stronger version of Theorem 2. To formulate it, let us define for $T \subset (\mathbb{R}^n)^{d-1}$,

$$F_A(T) := \mathbf{E} \sup_{\mathbf{x} \in T}\left(\sum_{\mathbf{i}} a_{\mathbf{i}} \prod_{k=1}^{d-1} x_{i_k}^k g_{i_d}\right).$$

THEOREM 3. *For any $T \subset (B_2^n)^{d-1}$ and $p \geq 1$,*

(17) $$F_A(T) \leq C(d)\left(\sqrt{p}\Delta_A(T) + \sum_{k=1}^{d-1} p^{(1-k)/2} s_{d-k}(A)\right).$$

Let us observe that $\mathbf{E}\|(\sum_{i_d} a_{\mathbf{i}} g_{i_d})\|_{\{1\}\ldots\{d-1\}} = F_A((B_2^n)^{d-1})$ and, therefore, Theorem 3 implies Theorem 2 since by (14), $\Delta_A((B_2^n)^{d-1}) \leq 2\|A\|_{\{1\}\ldots\{d\}}$.

We will prove (17) by induction on $d$, but first we will show several consequences of the theorem. In the next three lemmas, we shall assume that Theorem 3 (and thus also Theorem 2) holds for all matrices of order smaller than $d$.

For $\mathbf{x}, \mathbf{y} \in (\mathbb{R}^n)^{d-1}$, we set

$$\widetilde{\alpha}_A(\mathbf{x}) = \sum_{1 \leq j \neq k \leq d-1} \left\|\sum_{i_j} a_{\mathbf{i}} x_{i_j}^j\right\|_{\{k,d\}\{\{l\} : 1 \leq l \leq d-1, l \neq k,j\}},$$

$$\widetilde{\rho}_A(\mathbf{x}, \mathbf{y}) := \widetilde{\alpha}_A(\mathbf{x} - \mathbf{y}) = \sum_{1 \leq j \neq k \leq d-1} \left\|\sum_{i_j} a_{\mathbf{i}}(x_{i_j}^j - y_{i_j}^j)\right\|_{\{k,d\}\{\{l\} : 1 \leq l \leq d-1, l \neq k,j\}}$$

and for $T \subset (\mathbb{R}^n)^{d-1}$,

$$\widetilde{\alpha}_A(T) := \sup\{\widetilde{\alpha}_A(\mathbf{x}) : \mathbf{x} \in T\}.$$

LEMMA 6. *For any $p \geq 1$ and $l \geq 0$,*

$$N\left((B_2^n)^{d-1}, \widetilde{\rho}_A, 2^{-l}\sum_{k=1}^{d-1}p^{(1-k)/2}s_{d-k}(A)\right) \leq \exp(C(d)2^{2l}p).$$

PROOF. Note that $\widetilde{\alpha}_A$ is a norm on $(\mathbb{R}^n)^{d-1} = \mathbb{R}^{(d-1)n}$ and that

$$\mathbf{E}\widetilde{\alpha}(G_{(d-1)n}) = \sum_{1\leq j\neq k\leq d-1} \mathbf{E}\left\|\sum_{i_j} a_{\mathbf{i}}g_{i_j}\right\|_{\{k,d\}\{\{l\}:1\leq l\leq d-1, l\neq k,j\}}.$$

Let us fix $1 \leq j \neq k \leq d-1$ and observe that

$$\left\|\sum_{i_j} a_{\mathbf{i}}g_{i_j}\right\|_{\{k,d\}\{\{l\}:1\leq l\leq d-1, l\neq k,j\}} = \left\|\sum_{i_{d-1}} b_{i_1,\ldots,i_{d-1}}g_{i_{d-1}}\right\|_{\{1\}\ldots\{d-2\}}$$

for an appropriately chosen matrix $B = (b_{i_1,\ldots,i_{d-1}})$ (we treat a pair of indices $k, d$ as a single index and renumerate indices in such a way that $j$ would become $d-1$). Moreover, for any $1 \leq l \leq d-1$,

$$\sum_{(I_1,\ldots,I_l)\in S(l,d-1)} \|B\|_{I_1,\ldots,I_l} = \sum_{\substack{(I_1,\ldots,I_l)\in S(l,d)\\ \{k,d\}\in I_1}} \|A\|_{I_1,\ldots,I_l} \leq s_l(A).$$

Thus, by (5) (applied to the matrix $B$ of order $d-1$),

$$\mathbf{E}\left\|\sum_{i_j} a_{\mathbf{i}}g_{i_j}\right\|_{\{k,d\}\{\{l\}:1\leq l\leq d-1, l\neq k,j\}}$$

$$= \mathbf{E}\left\|\sum_{i_{d-1}} b_{\mathbf{i}}g_{i_{d-1}}\right\|_{\{1\}\ldots\{d-2\}}$$

$$\leq C(d)\sum_{l=1}^{d-1}p^{(2-d+l)/2}\sum_{(I_1,\ldots,I_l)\in S(l,d-1)}\|B\|_{I_1,\ldots,I_l}$$

$$\leq C(d)\sum_{s=1}^{d-1}p^{(2-s)/2}s_{d-s}(A).$$

Hence, by Corollary 2 (with $d=1$), we have for $t \in (0,1]$,

$$N\left((B_2^n)^{d-1}, \widetilde{\rho}_A, C(d)t\sum_{k=1}^{d-1}p^{(2-k)/2}s_{d-k}(A)\right) \leq \exp(Ct^{-2})$$

and it suffices to make the substitution $t = (C(d)2^l\sqrt{p})^{-1}$. □

LEMMA 7. *Suppose that $d \geq 3$, $\mathbf{y} \in (B_2^n)^{d-1}$ and $T \subset (B_2^n)^{d-1}$. Then for any $p \geq 1$ and $l \geq 0$, we can find a decomposition*

$$T = \bigcup_{j=1}^{N} T_j, \qquad N \leq \exp(C(d)2^{2l}p)$$

*such that for each $j \leq N$,*

(18)
$$F_A(\mathbf{y} + T_j) \leq F_A(T_j) + C(d)\left(\widetilde{\alpha}_A(\mathbf{y}) + \widetilde{\alpha}_A(T) + 2^{-l}\sum_{k=2}^{d-1} p^{(1-k)/2} s_{d-k}(A)\right)$$

*and*

(19)
$$\Delta_A(T_j) \leq 2^{-l} p^{-1/2} \widetilde{\alpha}_A(T) + 2^{-2l} \sum_{k=2}^{d-1} p^{-k/2} s_{d-k}(A).$$

PROOF. For $I \subsetneq \{1, \ldots, d-1\}$, $\mathbf{x}, \tilde{\mathbf{x}} \in (\mathbb{R}^n)^{d-1}$ and $S \subset (\mathbb{R}^n)^{d-1}$, let us define

$$\rho_A^{\mathbf{y},I}(\mathbf{x}, \tilde{\mathbf{x}}) := \left(\sum_{i_d}\left(\sum_{i_1,\ldots,i_{d-1}} a_{\mathbf{i}} \prod_{k \in I} y_{i_k}^k \left(\prod_{j \leq d-1, j \notin I} x_{i_j}^j - \prod_{j \leq d-1, j \notin I} \tilde{x}_{i_j}^j\right)\right)^2\right)^{1/2},$$

$$\Delta_A^{\mathbf{y},I}(S) := \sup\{\rho_A^{\mathbf{y},I}(\mathbf{x}, \tilde{\mathbf{x}}) : \mathbf{x}, \tilde{\mathbf{x}} \in S\}$$

and

$$F_A^{\mathbf{y},I}(S) := \mathbf{E}\sup_{\mathbf{x} \in S} \sum_{\mathbf{i}} a_{\mathbf{i}} \prod_{k \in I} y_{i_k}^k \prod_{j \leq d-1, j \notin I} x_{i_j}^j g_{i_d}.$$

Note that if $I \neq \emptyset$, then (17) applied to the matrix

$$A(\mathbf{y}, I) := \left(\sum_{i_I} a_{\mathbf{i}} \prod_{k \in I} y_{i_k}^k\right)$$

of order $d - \#I < d$ gives for any $S \subset (B_2^n)^{d-1}$ and $q \geq 1$,

$$F_A^{\mathbf{y},I}(S) \leq C(d - \#I)\left(q^{1/2}\Delta_A^{\mathbf{y},I}(S) + \sum_{k=1}^{d-\#I-1} q^{(1-k)/2} s_{d-\#I-k}(A(\mathbf{y},I))\right).$$

But, $s_{d-\#I-k}(A(\mathbf{y},I)) \leq s_{d-k}(A)$ for $k \geq 2$ and $s_{d-\#I-1}(A(\mathbf{y},I)) \leq \widetilde{\alpha}_A(\mathbf{y})$, hence,

(20) $$F_A^{\mathbf{y},I}(S) \leq C(d)\left(q^{1/2}\Delta_A^{\mathbf{y},I}(S) + \widetilde{\alpha}_A(\mathbf{y}) + \sum_{k=2}^{d-1} q^{(1-k)/2} s_{d-k}(A)\right).$$

Since $\mathbf{E}\sum_{\mathbf{i}} a_{\mathbf{i}} \prod_{k \leq d-1} y_{i_k}^k g_{i_d} = 0$, we get

$$(21) \qquad F_A(y+S) - F_A(S) \leq \sum_{\varnothing \neq I \subsetneq \{1,\ldots,d-1\}} F_A^{\mathbf{y},I}(S).$$

Observe also that for any $I \subset \{1,\ldots,d-1\}$, $0 \leq \#I \leq d-3$, we have

$$W_1^T(A(\mathbf{y},I)) \leq \sup\{\widetilde{\alpha}_A(\mathbf{x}) : \mathbf{x} \in T\} = \widetilde{\alpha}_A(T).$$

Thus, we may apply $2^{d-1} - d$ times (16) with $t = 2^{-l}p^{-1/2}$ and find a decomposition $T = \bigcup_{j=1}^N T_j$, $N \leq \exp(C(d)2^{2l}p)$, such that for each $j$ and $I \subset \{1,\ldots,d-1\}$, $0 \leq \#I \leq d-3$,

$$(22) \qquad \Delta_A^{\mathbf{y},I}(T_j) \leq 2^{-l}p^{-1/2}\widetilde{\alpha}_A(T) + 2^{-2l}\sum_{k=2}^{d-1} p^{-k/2}s_{d-k}(A).$$

Moreover, if $I \subset \{1,\ldots,d-1\}$ with $\#I = d-2$, then $A(\mathbf{y},I)$ is a matrix of order 2 and for $S \subset (B_2^n)^{d-1}$,

$$(23) \qquad F_A^{\mathbf{y},I}(S) \leq \|A(\mathbf{y},I)\|_{\mathrm{HS}} \leq \widetilde{\alpha}_A(\mathbf{y}).$$

Estimate (22) reduces to (19) for $I = \varnothing$ and (18) follows by (20) with $q = 2^{2l}p$ and (21)–(23). $\square$

LEMMA 8. *Suppose that $S$ is a finite subset of $(B_2^n)^{d-1}$, with $\#S \geq 2$, such that $S - S \subset (B_2^n)^{d-1}$. Then there exist finite sets $S_i \subset (B_2^n)^{d-1}$ and $\mathbf{y}_i \in S$, $i = 1,\ldots,N$, such that:*

(i) $2 \leq N \leq \exp(C(d)2^{2l}p)$,
(ii) $S = \bigcup_{i=1}^N (\mathbf{y}_i + S_i)$, $S_i - S_i \subset S - S$, $\#S_i \leq \#S - 1$,
(iii) $\Delta_A(S_i) \leq 2^{-2l}\sum_{k=1}^{d-1} p^{-k/2}s_{d-k}(A)$,
(iv) $\widetilde{\alpha}_A(S_i) \leq 2^{-l}\sum_{k=1}^{d-1} p^{(1-k)/2}s_{d-k}(A)$

*and*

(v) $F_A(\mathbf{y}_i + S_i) \leq F_A(S_i) + C(d)(\widetilde{\alpha}_A(S) + 2^{-l}\sum_{k=1}^{d-1} p^{(1-k)/2}s_{d-k}(A))$.

PROOF. By Lemma 6, we get

$$S = \bigcup_{i=1}^{N_1}(\mathbf{y}_i + T_i), \qquad N_1 \leq \exp(C(d)2^{2l}p),$$

$\mathbf{y}_i \in S$, $0 \in T_i$ and

$$\widetilde{\alpha}_A(T_i) \leq 2^{-l}\sum_{k=1}^{d-1} p^{(1-k)/2}s_{d-k}(A).$$

Note that $T_i \subset S - \mathbf{y}_i \subset S - S \subset (B_2^n)^{d-1}$. Hence, by Lemma 7 (with $l+1$ instead of $l$), we get

$$T_i = \bigcup_{j=1}^{N_2} T_{i,j}, \qquad N_2 \leq \exp(C(d)2^{2l}p),$$

where

$$F_A(\mathbf{y}_i + T_{i,j}) \leq F_A(T_{i,j}) + C(d)\left(\widetilde{\alpha}_A(\mathbf{y}_i) + \widetilde{\alpha}_A(T_i) + 2^{-l}\sum_{k=2}^{d-1} p^{(1-k)/2} s_{d-k}(A)\right)$$

$$\leq F_A(T_{i,j}) + C(d)\left(\widetilde{\alpha}_A(S) + 2^{-l}\sum_{k=2}^{d-1} p^{(1-k)/2} s_{d-k}(A)\right)$$

and

$$\Delta_A(T_{i,j}) \leq 2^{-l-1} p^{-1/2} \widetilde{\alpha}_A(T_i) + 2^{-2l-2}\sum_{k=2}^{d-1} p^{-k/2} s_{d-k}(A)$$

$$\leq 2^{-2l}\sum_{k=1}^{d-1} p^{-k/2} s_{d-k}(A).$$

Therefore,

$$S = \bigcup_{i,j}(\mathbf{y}_i + T_{i,j}).$$

We have $N = N_1 N_2 \leq \exp(C(d)2^{2l}p)$. Moreover, we may obviously assume that $N \geq 2$, and by making the sets $T_{i,j}$ disjoint, we may assume that $\#T_{i,j} \leq \#S - 1$. Obviously, $T_{i,j} - T_{i,j} \subset S - S$ and $\widetilde{\alpha}_A(T_{i,j}) \leq \widetilde{\alpha}_A(T_i) \leq 2^{-l}\sum_{k=1}^{d-1} p^{(1-k)/2} s_{d-k}(A)$. □

PROOF OF THEOREM 3. We proceed by induction on $d$. For $d=2$ and $A = (a_{ij})$, we have

$$F_A(T) \leq F_A(B_2^n) = \mathbf{E}\left(\sum_i \left(\sum_j a_{ij} g_j\right)^2\right)^{1/2} \leq \|A\|_{\mathrm{HS}} = s_1(A).$$

Suppose that $d \geq 3$ and that (17) holds for matrices of order smaller than $d$. Let us put $\Delta_0 := \Delta_A(T)$, $\widetilde{\Delta}_0 := \widetilde{\alpha}_A((B_2^n)^{d-1}) \leq C(d) s_{d-1}(A)$ and

$$\Delta_l := 2^{2-2l}\sum_{k=1}^{d-1} p^{-k/2} s_{d-k}(A),$$

$$\widetilde{\Delta}_l := 2^{1-l}\sum_{k=1}^{d-1} p^{(1-k)/2} s_{d-k}(A) \qquad \text{for } l \geq 1.$$

Suppose first that $T \subset \frac{1}{2}(B_2^n)^{d-1}$ and define

$$c_T(r,l) := \sup\{F_A(S) : S \subset (B_2^n)^{d-1}, S - S \subset T - T,$$
$$\#S \leq r, \Delta_A(S) \leq \Delta_l, \widetilde{\alpha}_A(S) \leq \widetilde{\Delta}_l\}.$$

Note that any subset $S \subset T$ satisfies $\Delta_A(S) \leq \Delta_0$ and $\widetilde{\alpha}_A(S) \leq \widetilde{\Delta}_0$, therefore,

(24) $$c_T(r,0) \geq \sup\{F_A(S) : S \subset T, \#S \leq r\}.$$

Obviously, $c_T(1,l) = 0$. We will now show that for $r \geq 2$,

(25)
$$c_T(r,l) \leq c_T(r-1, l+1)$$
$$+ C(d)\bigg(\widetilde{\Delta}_l + 2^l\sqrt{p}\Delta_l + 2^{-l}\sum_{k=1}^{d-1} p^{(1-k)/2} s_{d-k}(A)\bigg).$$

Indeed, let us take $S \subset (B_2^n)^{d-1}$ as in the definition of $c_T(r,l)$. Then by Lemma 8, we may find a decomposition $S = \bigcup_{i=1}^N (\mathbf{y}_i + S_i)$ satisfying (i)–(v). Hence, by Lemma 3 and (13), we have

(26)
$$F_A(S) \leq C\sqrt{\log N}\Delta_A(S) + \max_i F_A(\mathbf{y}_i + S_i)$$
$$\leq C(d)\bigg(\widetilde{\alpha}_A(S) + 2^l\sqrt{p}\Delta_l + 2^{-l}\sum_{k=1}^{d-1} p^{(1-k)/2} s_{d-k}(A)\bigg)$$
$$+ \max_i F_A(S_i).$$

We have $\Delta_A(S_i) \leq \Delta_{l+1}$, $\widetilde{\alpha}_A(S_i) \leq \widetilde{\Delta}_{l+1}$, $S_i - S_i \subset S - S \subset T - T$ and $\#S_i \leq \#S - 1 \leq r - 1$, thus $\max_i F_A(S_i) \leq c_T(r-1, l+1)$ and (26) yields (25).

By (25), we immediately obtain

$$c_T(r,0) \leq c_T(1, r-1) + C(d)\sum_{l=0}^{\infty}\bigg(\widetilde{\Delta}_l + 2^l\sqrt{p}\Delta_l + 2^{-l}\sum_{k=1}^{d-1} p^{(1-k)/2} s_{d-k}(A)\bigg)$$
$$\leq C(d)\bigg(\sqrt{p}\Delta_A(T) + \sum_{k=1}^{d-1} p^{(1-k)/2} s_{d-k}(A)\bigg).$$

For $T \subset \frac{1}{2}(B_2^n)^{d-1}$, we have by (24),

$$F_A(T) = \sup\{F_A(S) : S \subset T, \#S < \infty\} \leq \sup_r c_T(r,0)$$
$$\leq C(d)\bigg(\sqrt{p}\Delta_A(T) + \sum_{k=1}^{d-1} p^{(1-k)/2} s_{d-k}(A)\bigg).$$

Finally, if $T \subset (B_2^n)^{d-1}$, then $\frac{1}{2}T \subset \frac{1}{2}(B_2^n)^{d-1}$ and $\Delta_A(\frac{1}{2}T) = 2^{1-d}\Delta_A(T)$, hence,

$$F_A(T) = 2^{d-1}F_A(\tfrac{1}{2}T)$$
$$\leq C(d)\left(\sqrt{p}\Delta_A(T) + \sum_{k=1}^{d-1} p^{(1-k)/2} s_{d-k}(A)\right). \qquad \Box$$

PROOF OF THEOREM 1. First, we prove by induction on $d$ the estimate from below. For $d = 1$, it follows by (1). Suppose that the lower estimate holds for matrices of order smaller than $d$. Then by the induction assumption, we have for any matrix $B = (b_\mathbf{i})_{i_1,\ldots,i_{d-1}}$,

$$\left\|\sum_\mathbf{i} b_\mathbf{i} \prod_{j=1}^{d-1} g_{i_j}^{(j)}\right\|_p \geq C(d-1)^{-1}\sqrt{p}\left(\sum_\mathbf{i} b_\mathbf{i}^2\right)^{1/2} \geq C(d)^{-1}\left\|\sum_\mathbf{i} b_\mathbf{i} g_{i_1,\ldots,i_{d-1}}\right\|_p,$$

where $(g_{i_1,\ldots,i_{d-1}})$ is a sequence of i.i.d. $\mathcal{N}(0,1)$ r.v.'s independent of $(g_{i_d}^{(d)})$. Therefore, by (2),

(27)
$$\left\|\sum_\mathbf{i} a_\mathbf{i} \prod_{j=1}^{d} g_{i_j}^{(j)}\right\|_p \geq C(d)^{-1}\left\|\sum_\mathbf{i} a_\mathbf{i} g_{i_1,\ldots,i_{d-1}} g_{i_d}^{(d)}\right\|_p$$
$$\geq C(d)^{-1}\sqrt{p}\|A\|_{\{1,\ldots,d\}}.$$

Let $(I_1,\ldots,I_k) \in S(k,d)$ with $k \geq 2$ and $\sum_{i_{I_l}} (x_{I_l}^{(l)})^2 \leq 1$ for $l = 1,\ldots,k$. Then by the induction assumption applied twice [first conditionally on $(g_{i_j}^{(j)})_{j\in I_1}$], we have

$$\left\|\sum_\mathbf{i} a_\mathbf{i} \prod_{j=1}^{d} g_{i_j}^{(j)}\right\|_p \geq C(d - \#I_1)^{-1} p^{(k-1)/2}\left\|\sum_\mathbf{i} a_\mathbf{i} \prod_{j\in I_1} g_{i_j}^{(j)} \prod_{l=2}^{k} x_{I_l}^{(l)}\right\|_p$$
$$\geq (C(d - \#I_1)C(\#I_1))^{-1} p^{k/2} \sum_\mathbf{i} a_\mathbf{i} \prod_{l=1}^{k} x_{I_l}^{(l)}$$

and, hence,

(28)
$$\left\|\sum_\mathbf{i} a_\mathbf{i} \prod_{j=1}^{d} g_{i_j}^{(j)}\right\|_p \geq C(d)^{-1} p^{k/2}\|A\|_{I_1,\ldots,I_k}.$$

Inequalities (27) and (28) imply the lower part of estimate (4).

Now (again by induction on $d$), we prove the estimate from above. For $d \leq 2$, the estimate follows by (1) and (2). Suppose that $d \geq 3$ and the estimate holds for chaoses of order smaller than $d$. By the induction assumption, we have

$$
(29) \quad \left\|\sum_{\mathbf{i}} a_{\mathbf{i}} \prod_{j=1}^{d} g_{i_j}^{(j)}\right\|_p \\
\leq C(d-1) \sum_{k=1}^{d-1} p^{k/2} \sum_{(I_1,\ldots,I_k)\in S(k,d-1)} \left(\mathbf{E}\left\|\left(\sum_{i_d} a_{\mathbf{i}} g_{i_d}^{(d)}\right)\right\|_{I_1,\ldots,I_k}^p\right)^{1/p}.
$$

However, for $(I_1,\ldots,I_k) \in S(k,d-1)$, we have by (12),

$$
(30) \quad \left(\mathbf{E}\left\|\left(\sum_{i_d} a_{\mathbf{i}} g_{i_d}^{(d)}\right)\right\|_{I_1,\ldots,I_k}^p\right)^{1/p} \\
\leq C\sqrt{p}\|A\|_{I_1,\ldots,I_k\{d\}} + \mathbf{E}\left\|\left(\sum_{i_d} a_{\mathbf{i}} g_{i_d}^{(d)}\right)\right\|_{I_1,\ldots,I_k}.
$$

Theorem 2 gives

$$
(31) \quad \mathbf{E}\left\|\left(\sum_{i_d} a_{\mathbf{i}} g_{i_d}^{(d)}\right)\right\|_{I_1,\ldots,I_k} \leq C(k+1) p^{-k/2} m_p(A).
$$

Inequalities (29)–(31) then yield the upper estimate in (4). □

PROOF OF COROLLARY 1. Let $S := \sum_{\mathbf{i}} a_{\mathbf{i}} \prod_{j=1}^{d} g_{i_j}^{(j)}$. By Chebyshev's inequality and (4), one gets for $p \geq 2$,

$$
(32) \quad \mathbf{P}(|S| \geq eC(d) m_p(A)) \leq \mathbf{P}(|S| \geq e\|S\|_p) \leq e^{-p}.
$$

Since $\|S\|_{2p} \leq C_1(d)\|S\|_p$ (cf. [9], Section 3.2, or use (4) and $m_{2p}(A) \leq 2^{d/2} m_p(A)$), we get by the Paley–Zygmund inequality for $q \geq 2$,

$$
\mathbf{P}(S \geq 2^{-1}\|S\|_q) = \mathbf{P}(|S|^q \geq 2^{-q} \mathbf{E}|S|^q) \\
\geq (1 - 2^{-q})^2 \frac{(\mathbf{E}|S|^q)^2}{\mathbf{E}|S|^{2q}} \geq (2C_1(d))^{-2q}.
$$

Fix $p > 0$ and take $q := p/(2\ln(2C_1(d)))$. Then $q \geq 2$ for $p \geq p_0(d)$ and $\|S\|_q \geq C(d)^{-1} m_q(A) \geq C(d)^{-1}(\max(2\ln(2C_1(d)),1))^{-d/2} m_p(A) = 2C_2(d)^{-1} m_p(A)$. Thus, for $p \geq p_0(d)$,

$$
\mathbf{P}(|S| \geq C_2(d)^{-1} m_p(A)) \geq \mathbf{P}(S \geq 2^{-1}\|S\|_q) \geq (2C_1(d))^{-2q} = e^{-p}
$$

and, therefore, for any $p > 0$,

$$
(33) \quad \mathbf{P}(|S| \geq C_2(d)^{-1} m_p(A)) \geq \min(e^{-p_0(d)}, e^{-p}).
$$

Finally, note that if $m_p(A) = s \geq m_2(A)$, then $p$ is comparable (with constants depending only on $d$) with

$$\min\{(s/\|A\|_{I_1,\ldots,I_k})^{2/k} : 1 \leq k \leq d, (I_1,\ldots,I_k) \in S(k,d)\}$$

and, therefore, Corollary 1 follows by (32) and (33). $\square$

**Acknowledgment.** The author would like to thank the referee for a careful reading of the manuscript and valuable remarks.

INSTITUTE OF MATHEMATICS
WARSAW UNIVERSITY
BANACHA 2
02-097 WARSZAWA
POLAND
E-MAIL: rlatala@mimuw.edu.pl